\documentclass{gtart}

\def\ifplaintex{\expandafter\ifx\csname documentclass\endcsname\relax}

\ifplaintex 
\else
\fi

\expandafter\ifx\csname beginpicture\endcsname\relax
\expandafter\ifx\csname documentclass\endcsname\relax
\input pictex \else
\input prepictex \input pictex \input postpictex \fi\fi

\def\gt{{\mathsurround=0pt\it $\cal G\mskip-2mu$eometry \&\ 
$\cal T\!\!$opology}}        

\def\gtp{{\mathsurround=0pt\it $\cal G\mskip-2mu$eometry \&\ 
$\cal T\!\!$opology $\cal P\!$ublications}}  

\def\lognumber#1{\def\thelognumber{#1}}
\def\volumenumber#1{\def\thevolumenumber{#1}}
\def\papernumber#1{\def\thepapernumber{#1}}
\def\volumeyear#1{\def\thevolumeyear{#1}}

\def\pagenumbers#1#2{\def\startpage{#1}\def\finishpage{#2}}
\def\published#1{\def\publishdate{#1}}
\def\proposed#1{\def\theproposer{#1}}
\def\seconded#1{\def\theseconders{#1}}
\def\received#1{\def\receiveddate{#1}}
\def\revised#1{\def\reviseddate{#1}}
\def\accepted#1{\def\accepteddate{#1}}

\long\def\asciiabstract#1{\long\def\theasciiabstract{#1}}

\def\shorttitle#1{\def\theshorttitle{#1}}

\let\\\par\let\thelognumber\relax
\let\thevolumenumber\relax\let\thepapernumber\relax
\let\thevolumeyear\relax\let\thesamplenumber\relax\let\startpage\relax
\let\finishpage\relax\let\publishdate\relax\let\receiveddate\relax
\let\reviseddate\relax\let\accepteddate\relax\let\theasciititle\relax
\let\theasciiauthors\relax
\let\theasciiabstract\relax
\let\theasciiemail\relax\let\theshortauthors\relax\let\theshorttitle\relax

\long\def\maketitlep{   

\count0=\startpage

\gt\hfill      
\beginpicture
\setcoordinatesystem units <0.33truein, 0.33truein> point at 2.2 0.9
\setplotsymbol ({$\cal G$})
\plotsymbolspacing=9truept
\circulararc 315 degrees from 0 1 center at 0 0
\setplotsymbol ({$\cal T$})
\circulararc 315 degrees from 1 -1 center at 1 0
\endpicture
\break
{\small\ifx\thesamplenumber\relax 
Volume \else Sample
\fi\thevolumenumber\ (\thevolumeyear)
\startpage--\finishpage\nl
Published: \publishdate}
\vglue 0.5truein plus 0.4fil minus 0.1truein

{\parskip=0pt\leftskip 0pt plus 1fil\def\\{\par\smallskip}{\ifplaintex\large
\else\Large\fi\bf\thetitle}\par\medskip}   

\vglue 0pt plus 0.1fil 

{\parskip=0pt\leftskip 0pt plus 1fil\def\\{\par}{\sc\theauthors}
\par\medskip}

\vglue 0pt plus 0.1fil 

{\small\parskip=0pt\let\newline\\
{\leftskip 0pt plus 1fil\def\\{\par}{\sl\theaddress}\par}
\expandafter\ifx\theemail\relax    
\relax\else\vglue 5pt plus 0.02fil minus 2pt\def\\{\stdspace{\rm 
and}\stdspace} 
\cl{Email:\stdspace\tt\theemail}\fi
\ifx\theurl\relax                  
\relax\else\vglue 5pt plus 0.02fil minus 2pt\def\\{\stdspace{\rm 
and}\stdspace}
\cl{URL:\stdspace\tt\theurl}\fi\par}

\vglue 7pt plus 0.3fil minus 3pt

{\bf Abstract}
\vglue 5pt plus 0.1fil minus 2pt

\theabstract

\vglue 7pt plus 0.3fil minus 3pt

{\bf AMS Classification numbers}\quad Primary:\quad \theprimaryclass

Secondary:\quad \thesecondaryclass

\vglue 5pt plus 0.3fil minus 2pt

{\bf Keywords}\quad \thekeywords

\vglue 10pt plus 0.5fil minus 5pt

{\small  Proposed: \theproposer\hfill Received: \receiveddate\nl
Seconded: \theseconders\hfill 
\ifx\reviseddate\relax                         
Accepted: \accepteddate                        
\else
Revised: \reviseddate                          
\fi}
\eject
}       

\let\maketitlepage\maketitlep
\let\maketitle\maketitlepage

\font\phead=cmsl9 scaled 950
\font=cmsl9 scaled 1050
\font\pnum=cmbx10 scaled 913
\font\lnum=cmbx10 
\font\pfoot=cmsl9 scaled 950
\font=cmsl9 scaled 1050
\ifplaintex
\headline{\vbox to 0pt{\vskip -4.5mm\line{\small\phead\ifnum
\count0=\startpage ISSN 1364-0380 (on line)
1465-3060 (printed) \hfill {\pnum\folio}\else\ifodd\count0\def\\{ }%
\ifx\theshorttitle\relax\thetitle\else\theshorttitle\fi\hfill{\pnum\folio}
\else\def\\{ and }{\pnum\folio}\hfill\ifx\theshortauthors\relax\theauthors
\else\theshortauthors\fi\fi\fi}\vss}}
\footline{\vbox to 0pt{\vglue 0mm\line{\small\pfoot\ifnum\count0=\startpage
\copyright\ \gtp\hfill\else
\gt, Volume \thevolumenumber\ (\thevolumeyear)\hfill\fi}\vss
}}
\else
\makeatletter
\def\@oddhead{{\small\lhead\ifnum\count0=\startpage ISSN 1364-0380 (on line)
1465-3060 (printed) \hfill {\lnum\number\count0}\else\ifodd\count0
\def\\{ }\ifx\theshorttitle\relax \thetitle \else\theshorttitle\fi\hfill
{\lnum\number\count0}\else\def\\{ and }{\lnum\number\count0}
\hfill\ifx\theshortauthors\relax 
\theauthors\else\theshortauthors\fi\fi\fi}}\def\@evenhead{\@oddhead}
\def\@oddfoot{\small\lfoot\ifnum\count0=\startpage\copyright\ \gtp\hfill\else
\gt, Volume \thevolumenumber\ (\thevolumeyear)\hfill\fi}
\def\@evenfoot{\@oddfoot}
\makeatother
\fi

\newwrite\gtoutfile
\long\gdef\makeheadfile{  
{\def\\{, }\def\s{ }
\immediate\openout\gtoutfile head.xxx
\immediate\write\gtoutfile{To: math@arxiv.org}
\immediate\write\gtoutfile{Subject: put or rep NNNNN:pppp}
\immediate\write\gtoutfile{--text follows this line--}
\immediate\write\gtoutfile{Proxy-for: \ifx\theasciiauthors\relax
\theauthors\else\theasciiauthors\fi\s<\ifx\theasciiemail\relax\theemail\else\theasciiemail\fi>}
\immediate\write\gtoutfile{\noexpand\\}
\immediate\write\gtoutfile{Authors: \ifx\theasciiauthors\relax
\theauthors\else\theasciiauthors\fi}
\immediate\write\gtoutfile{Title: \ifx\theasciititle\relax
\thetitle\else\theasciititle\fi}
\immediate\write\gtoutfile{Subj-class: GT or SG or MG etc}
\immediate\write\gtoutfile{MSC-class: \theprimaryclass\ifx\thesecondaryclass\relax\else, \thesecondaryclass\fi}
\immediate\write\gtoutfile{Journal-ref: Geom. Topol. \thevolumenumber
(\thevolumeyear) \startpage-\finishpage}
\immediate\write\gtoutfile{Comments: Published by Geometry and Topology at}
\immediate\write\gtoutfile{\s\s http://www.maths.warwick.ac.uk/gt/GTVol\thevolumenumber/paper\thepapernumber.abs.html}
\immediate\write\gtoutfile{\noexpand\\}
\immediate\write\gtoutfile{}
\ifx\theasciiabstract\relax
\immediate\write\gtoutfile{\theabstract}\else
\immediate\write\gtoutfile{\theasciiabstract}\fi
\immediate\write\gtoutfile{}
\immediate\write\gtoutfile{\noexpand\\}
\immediate\write\gtoutfile{}
\immediate\closeout\gtoutfile}}  

\def\maketitlepage{\maketitlep\makeheadfile}
\let\maketitle\maketitlepage

\def\ifplaintex{\expandafter\ifx\csname documentclass\endcsname\relax}

\ifplaintex 
\else
\fi

\expandafter\ifx\csname beginpicture\endcsname\relax
\expandafter\ifx\csname documentclass\endcsname\relax
\input pictex \else
\input prepictex \input pictex \input postpictex \fi\fi

\def\gt{{\mathsurround=0pt\it $\cal G\mskip-2mu$eometry \&\ 
$\cal T\!\!$opology}}        

\def\gtp{{\mathsurround=0pt\it $\cal G\mskip-2mu$eometry \&\ 
$\cal T\!\!$opology $\cal P\!$ublications}}  

\def\lognumber#1{\def\thelognumber{#1}}
\def\volumenumber#1{\def\thevolumenumber{#1}}
\def\papernumber#1{\def\thepapernumber{#1}}
\def\volumeyear#1{\def\thevolumeyear{#1}}

\def\pagenumbers#1#2{\def\startpage{#1}\def\finishpage{#2}}
\def\published#1{\def\publishdate{#1}}
\def\proposed#1{\def\theproposer{#1}}
\def\seconded#1{\def\theseconders{#1}}
\def\received#1{\def\receiveddate{#1}}
\def\revised#1{\def\reviseddate{#1}}
\def\accepted#1{\def\accepteddate{#1}}

\long\def\asciiabstract#1{\long\def\theasciiabstract{#1}}

\def\shorttitle#1{\def\theshorttitle{#1}}

\let\\\par\let\thelognumber\relax
\let\thevolumenumber\relax\let\thepapernumber\relax
\let\thevolumeyear\relax\let\thesamplenumber\relax\let\startpage\relax
\let\finishpage\relax\let\publishdate\relax\let\receiveddate\relax
\let\reviseddate\relax\let\accepteddate\relax\let\theasciititle\relax
\let\theasciiauthors\relax
\let\theasciiabstract\relax
\let\theasciiemail\relax\let\theshortauthors\relax\let\theshorttitle\relax

\long\def\maketitlep{   

\count0=\startpage

\gt\hfill      
\beginpicture
\setcoordinatesystem units <0.33truein, 0.33truein> point at 2.2 0.9
\setplotsymbol ({$\cal G$})
\plotsymbolspacing=9truept
\circulararc 315 degrees from 0 1 center at 0 0
\setplotsymbol ({$\cal T$})
\circulararc 315 degrees from 1 -1 center at 1 0
\endpicture
\break
{\small\ifx\thesamplenumber\relax 
Volume \else Sample
\fi\thevolumenumber\ (\thevolumeyear)
\startpage--\finishpage\nl
Published: \publishdate}
\vglue 0.5truein plus 0.4fil minus 0.1truein

{\parskip=0pt\leftskip 0pt plus 1fil\def\\{\par\smallskip}{\ifplaintex\large
\else\Large\fi\bf\thetitle}\par\medskip}   

\vglue 0pt plus 0.1fil 

{\parskip=0pt\leftskip 0pt plus 1fil\def\\{\par}{\sc\theauthors}
\par\medskip}

\vglue 0pt plus 0.1fil 

{\small\parskip=0pt\let\newline\\
{\leftskip 0pt plus 1fil\def\\{\par}{\sl\theaddress}\par}
\expandafter\ifx\theemail\relax    
\relax\else\vglue 5pt plus 0.02fil minus 2pt\def\\{\stdspace{\rm 
and}\stdspace} 
\cl{Email:\stdspace\tt\theemail}\fi
\ifx\theurl\relax                  
\relax\else\vglue 5pt plus 0.02fil minus 2pt\def\\{\stdspace{\rm 
and}\stdspace}
\cl{URL:\stdspace\tt\theurl}\fi\par}

\vglue 7pt plus 0.3fil minus 3pt

{\bf Abstract}
\vglue 5pt plus 0.1fil minus 2pt

\theabstract

\vglue 7pt plus 0.3fil minus 3pt

{\bf AMS Classification numbers}\quad Primary:\quad \theprimaryclass

Secondary:\quad \thesecondaryclass

\vglue 5pt plus 0.3fil minus 2pt

{\bf Keywords}\quad \thekeywords

\vglue 10pt plus 0.5fil minus 5pt

{\small  Proposed: \theproposer\hfill Received: \receiveddate\nl
Seconded: \theseconders\hfill 
\ifx\reviseddate\relax                         
Accepted: \accepteddate                        
\else
Revised: \reviseddate                          
\fi}
\eject
}       

\let\maketitlepage\maketitlep
\let\maketitle\maketitlepage

\font\phead=cmsl9 scaled 950
\font=cmsl9 scaled 1050
\font\pnum=cmbx10 scaled 913
\font\lnum=cmbx10 
\font\pfoot=cmsl9 scaled 950
\font=cmsl9 scaled 1050
\ifplaintex
\headline{\vbox to 0pt{\vskip -4.5mm\line{\small\phead\ifnum
\count0=\startpage ISSN 1364-0380 (on line)
1465-3060 (printed) \hfill {\pnum\folio}\else\ifodd\count0\def\\{ }%
\ifx\theshorttitle\relax\thetitle\else\theshorttitle\fi\hfill{\pnum\folio}
\else\def\\{ and }{\pnum\folio}\hfill\ifx\theshortauthors\relax\theauthors
\else\theshortauthors\fi\fi\fi}\vss}}
\footline{\vbox to 0pt{\vglue 0mm\line{\small\pfoot\ifnum\count0=\startpage
\copyright\ \gtp\hfill\else
\gt, Volume \thevolumenumber\ (\thevolumeyear)\hfill\fi}\vss
}}
\else
\makeatletter
\def\@oddhead{{\small\lhead\ifnum\count0=\startpage ISSN 1364-0380 (on line)
1465-3060 (printed) \hfill {\lnum\number\count0}\else\ifodd\count0
\def\\{ }\ifx\theshorttitle\relax \thetitle \else\theshorttitle\fi\hfill
{\lnum\number\count0}\else\def\\{ and }{\lnum\number\count0}
\hfill\ifx\theshortauthors\relax 
\theauthors\else\theshortauthors\fi\fi\fi}}\def\@evenhead{\@oddhead}
\def\@oddfoot{\small\lfoot\ifnum\count0=\startpage\copyright\ \gtp\hfill\else
\gt, Volume \thevolumenumber\ (\thevolumeyear)\hfill\fi}
\def\@evenfoot{\@oddfoot}
\makeatother
\fi

\newwrite\gtoutfile
\long\gdef\makeheadfile{  
{\def\\{, }\def\s{ }
\immediate\openout\gtoutfile head.xxx
\immediate\write\gtoutfile{To: math@arxiv.org}
\immediate\write\gtoutfile{Subject: put or rep NNNNN:pppp}
\immediate\write\gtoutfile{--text follows this line--}
\immediate\write\gtoutfile{Proxy-for: \ifx\theasciiauthors\relax
\theauthors\else\theasciiauthors\fi\s<\ifx\theasciiemail\relax\theemail\else\theasciiemail\fi>}
\immediate\write\gtoutfile{\noexpand\\}
\immediate\write\gtoutfile{Authors: \ifx\theasciiauthors\relax
\theauthors\else\theasciiauthors\fi}
\immediate\write\gtoutfile{Title: \ifx\theasciititle\relax
\thetitle\else\theasciititle\fi}
\immediate\write\gtoutfile{Subj-class: GT or SG or MG etc}
\immediate\write\gtoutfile{MSC-class: \theprimaryclass\ifx\thesecondaryclass\relax\else, \thesecondaryclass\fi}
\immediate\write\gtoutfile{Journal-ref: Geom. Topol. \thevolumenumber
(\thevolumeyear) \startpage-\finishpage}
\immediate\write\gtoutfile{Comments: Published by Geometry and Topology at}
\immediate\write\gtoutfile{\s\s http://www.maths.warwick.ac.uk/gt/GTVol\thevolumenumber/paper\thepapernumber.abs.html}
\immediate\write\gtoutfile{\noexpand\\}
\immediate\write\gtoutfile{}
\ifx\theasciiabstract\relax
\immediate\write\gtoutfile{\theabstract}\else
\immediate\write\gtoutfile{\theasciiabstract}\fi
\immediate\write\gtoutfile{}
\immediate\write\gtoutfile{\noexpand\\}
\immediate\write\gtoutfile{}
\immediate\closeout\gtoutfile}}  

\def\maketitlepage{\maketitlep\makeheadfile}
\let\maketitle\maketitlepage

\lognumber{215}
\volumenumber{5}\papernumber{27}\volumeyear{2001}
\pagenumbers{885}{893}
\received{26 October 2001}
\accepted{26 November 2001}
\proposed{Walter Neumann}
\seconded{Ralph Cohen, Steven Ferry}
\revised{10 November 2001}
\published{26 November 2001}

\usepackage{amssymb}
\usepackage{amsmath}
\newtheorem{thm}{Theorem}
\newcommand{\av}{\operatorname{av}}
\begin{document}

\title{A proof of Atiyah's conjecture on configurations\\of 
four points in Euclidean three-space}
\authors{Michael Eastwood\\Paul Norbury}
\shorttitle{Atiyah's Conjecture}
\address{Pure Mathematics Department\\Adelaide University\\
South Australia 5005}
\email{meastwoo@maths.adelaide.edu.au, pnorbury@maths.adelaide.edu.au}
\begin{abstract}
From any configuration of finitely many points in Euclidean three-space,
Atiyah constructed a determinant and conjectured that it was always non-zero.
In this article we prove the conjecture for the case of four points.
\end{abstract}
\asciiabstract{ 
From any configuration of finitely many points in Euclidean
three-space, Atiyah constructed a determinant and conjectured that it
was always non-zero.  Atiyah and Sutcliffe (hep-th/0105179) amass a
great deal of evidence it its favour.  In this article we prove the
conjecture for the case of four points.}

\maketitle
\renewcommand{\thefootnote}{$\sharp$}

Consider $n$ distinct points in Euclidean three-space.  Fixing
attention on one of these points, the others give rise to $n-1$ points
on its sphere of vision.  Thinking of this as the Riemann sphere gives
a monic polynomial of degree $\leq n-1$, having as its zeroes the
points not equal to the point chosen to be $\infty$.  We may regard
its coefficients as a complex $n$-vector (for degree $d<n-1$, its
first $n-1-d$ coefficients are deemed to be zero).  Repeating this
exercise for each of the $n$ points gives $n$ such vectors and hence
an $n\times n$ matrix.  In~\cite{a1,a2}, Atiyah conjectured that a
matrix constructed in this way cannot be singular.  In~\cite{as},
Atiyah and Sutcliffe amass a great deal of numerical evidence for this
conjecture and formulate a series of further conjectures based on the
geometry that their numerical studies apparently reveal.

In spite of overwhelming evidence in its favour, the basic conjecture,
as stated above, remains surprisingly resistant. The case $n=3$ is not
too hard: a geometric argument is given in \cite{a1} and an algebraic
one in~\cite{a2}.  In this article we establish the case $n=4$.

\section{Normalisation}
In describing Atiyah's conjecture above we used only the directions
defined by pairs of points amongst the $n$ points.  It turns out to be
far more natural to keep track of scale as well as direction, in
particular in order to see what happens if we rotate the sphere of
directions in~${\mathbb R}^3$, obtaining a different identification
with the Riemann sphere.  We will use the Hopf mapping
$${\mathbb C}^2\ni\left\lgroup\!\begin{array}c w_1\\ w_2
\end{array}\!\right\rgroup \stackrel{h}{\longmapsto}
\left\lgroup\!\begin{array}c (|w_1|^2-|w_2|^2)/2\\
w_1\bar{w}_2\end{array}\!\right\rgroup \in{\mathbb R}\times{\mathbb
C}\cong{\mathbb R}^3$$
as follows.  This mapping intertwines the action of ${\mathrm{SU}}(2)$
on ${\mathbb C}^2$ with the action of ${\mathrm{SO}}(3)$ on~${\mathbb
R}^3$ and descends to an isomorphism from ${\mathbb{CP}}_1$ to the
sphere of rays through the origin in ${\mathbb R}^3$.  Therefore, for
each point in ${\mathbb R}^3\setminus\{0\}$, we may choose a
corresponding point in ${\mathbb C}^2\setminus\{0\}$ defined up to
phase.  Their symmetric
tensor product lies in
$$\textstyle\bigodot^{n-1}{\mathbb C}^2\cong{\mathbb C}^n$$ and is
also well-defined up to phase.  We may regard this as normalising, up
to phase, the complex $n$-vectors appearing in our initial formulation
of Atiyah's conjecture.  If we now construct the columns of an
$n\times n$ matrix $M$ in this way, then $\det M$ is well-defined up
to phase and $|\det M|^2$ is invariant under Euclidean motions.

Following Atiyah and Sutcliffe~\cite{as}, we may normalize $\det M$ further.
Consider the mapping
$${\mathbb C}^2\ni\left\lgroup\!\begin{array}c w_1\\ w_2
\end{array}\!\right\rgroup
\stackrel{\sigma}{\longmapsto}\left\lgroup\!
\begin{array}c -\bar{w}_2\\ \bar{w}_1\end{array}\!\right\rgroup
\in{\mathbb C}^2,$$
observing that $h(\sigma(w))=-h(w)$ for all~$w\in{\mathbb C}^2$.
Also note that
\begin{equation}\label{phase}
\sigma(e^{i\theta}w)=e^{-i\theta}\sigma(w)\quad\mbox{and}\quad
\sigma(\sigma(w))=-w.
\end{equation}
Fix an ordering for our original $n$ points in ${\mathbb R}^3$.  Each
pair of these points contributes twice to $\det M$, once when the
later point is viewed from the earlier and once when this view is
reversed.  We mandate using $w$ and $\sigma(w)$, respectively, in
lifting to ${\mathbb C}^2$.  By virtue of~(\ref{phase}), both the
phase ambiguity $w\mapsto e^{i\theta}w$ and the ordering ambiguity
cancel from~$\det M$.  In conclusion, $\det M$ is invariant under
Euclidean motions.  It is easy to check that $\det M$ is replaced by
its complex conjugate under reflection.  In particular, if all points
lie in a plane then the determinant is real.  For further details
see~\cite{a2,as}.  We shall call~$\det M$, normalised in this way, the
{\em Atiyah determinant}.  In~\cite{as} a scale invariant
normalisation $D$ is used.  The two normalisations are related by
\[\det M=D\cdot\prod_{i>j}(2r_{ij})\] where $r_{ij}$ is the distance between
the $i^{\mathrm{th}}$ and $j^{\mathrm{th}}$ points.

We are free to use Euclidean motions to place points in convenient
locations.  Let us do this to verify the conjecture when $n=3$,
choosing the three points in ${\mathbb R}\times{\mathbb C}$ to be
$$\left\lgroup\!\begin{array}c0\\ 0\end{array}\!\right\rgroup\quad
\left\lgroup\!\begin{array}c0\\ a\end{array}\!\right\rgroup\quad
\left\lgroup\!\begin{array}c0\\ z\end{array}\!\right\rgroup$$
with $a$ real.
They form a triangle with side lengths
$a$, $b=|z|$, and $c=|a-z|$.
We may use the following
$$\begin{array}{ll}
\frac{1}{\sqrt a}\left\lgroup\!\begin{array}c a\\ a\end{array}\!\right\rgroup
\stackrel{h}{\mapsto}
\left\lgroup\!\begin{array}c0\\a\end{array}\!\right\rgroup&
\frac{1}{\sqrt b}\left\lgroup\!\begin{array}c b\\ \bar{z}\end{array}
\!\right\rgroup\stackrel{h}{\mapsto}
\left\lgroup\!\begin{array}c0\\ z\end{array}\!\right\rgroup\\ \\[-5pt]
\frac{1}{\sqrt a}\left\lgroup\!\begin{array}c -a\\ a\end{array}\!\right\rgroup
\stackrel{h}{\mapsto}
\left\lgroup\!\begin{array}c0\\-a\end{array}\!\right\rgroup&
\frac{1}{\sqrt c}\left\lgroup\!\begin{array}c c\\ \bar{z}-a
\end{array}\!\right\rgroup\stackrel{h}{\mapsto}
\left\lgroup\!\begin{array}c0\\ z-a\end{array}\!\right\rgroup\\ \\[-5pt]
\frac{1}{\sqrt b}\left\lgroup\!\begin{array}c -z\\ b\end{array}
\!\right\rgroup\stackrel{h}{\mapsto}
\left\lgroup\!\begin{array}c0\\ -z\end{array}\!\right\rgroup&
\frac{1}{\sqrt c}\left\lgroup\!\begin{array}c a-z\\
c\end{array}\!\right\rgroup
\stackrel{h}{\mapsto}
\left\lgroup\!\begin{array}c0\\a-z\end{array}\!\right\rgroup
\end{array}$$
in computing $\det M$.  We obtain
$$\begin{array}l
\displaystyle\frac{1}{abc}\left|\begin{array}{ccc}
ab&-ac&-z(a-z)\\
a\bar{z}+ab&-a(\bar{z}-a)+ac& -zc+b(a-z)\\
a\bar{z}&a(\bar{z}-a)&bc
\end{array}\right|\\ \\[-8pt]
\quad=a((z+\bar{z})(c-a-b)+2b(a+b+3c))\\[4pt]
\qquad=(a^2+b^2-c^2)(c-a-b)+2ab(a+b+3c)\\[4pt]
\quad\qquad=d_3(a,b,c)+8abc,\end{array}$$
where
\begin{equation}\label{deethree}d_3(a,b,c)=(a+b-c)(b+c-a)(c+a-b).\end{equation}
The triangle inequalities imply that $d_3(a,b,c)\geq 0$ with equality if and
only if the points lie on a line.  Therefore $\det M\geq 8abc>0$ and, in
particular, is non-zero, as required.

\section{The case $n=4$}
\begin{thm}
For any four points in ${\mathbb R}^3$, the Atiyah determinant is
non-zero.
\end{thm}
\begin{proof}
Choose the four points in ${\mathbb R}^3={\mathbb R}\times{\mathbb C}$
to be
$$\left\lgroup\!\begin{array}c0\\ z_1\end{array}\!\right\rgroup\quad
\left\lgroup\!\begin{array}c0\\ z_2\end{array}\!\right\rgroup\quad
\left\lgroup\!\begin{array}c0\\ z_3\end{array}\!\right\rgroup\quad
\left\lgroup\!\begin{array}cr\\ 0\end{array}\!\right\rgroup.$$ Put
$z_{ij}=z_i-z_j$ for $i>j$ and label the distances between points
by~$r_{ij}$.  We define $z_4=0$ so that $z_{4j}=-z_j$.  Thus,
$r_{ij}^2=|z_{ij}|^2$ for $i<4$ and $r_{4j}^2=r^2+|z_{4j}|^2$.

For $j<i<4$ the vector running from the $j^{\mathrm{th}}$ point to the
$i^{\mathrm{th}}$ point may be lifted to
\[w=\frac{1}{\sqrt{r_{ij}}}\left\lgroup\!\begin{array}c r_{ij}\\
\bar{z}_{ij} \end{array}\!\right\rgroup\mbox{, with }\sigma(w)=
\frac{1}{\sqrt{r_{ij}}}\left\lgroup\!\begin{array}c -z_{ij}\\
r_{ij} \end{array}\!\right\rgroup.\]
Similarly, if we put $R_{4j}=r_{4j}+r$ for $j<4$, then
\[w=\frac{1}{\sqrt{R_{4j}}}\left\lgroup\!\begin{array}cR_{4j}\\
\bar{z}_{4j}\end{array}\!\right\rgroup\mbox{ and }
\sigma(w)=\frac{1}{\sqrt{R_{4j}}}\left\lgroup\!\begin{array}c -z_{4j}\\R_{4j}
\end{array}\!\right\rgroup\]
lift to ${\mathbb C}^2$ the vectors in ${\mathbb R}^3$ joining the
$j^{\mathrm{th}}$ point to the $4^{\mathrm{th}}$ point and vice versa.
For each of the four points, the coefficients of the corresponding third
degree polynomial give the following four vectors:
\[v_1=\frac{1}{\sqrt{r_{21}r_{31}R_{41}}}\left\lgroup\!\begin{array}{c}
r_{21}r_{31}R_{41}\\
r_{21}r_{31}\bar{z}_{41}+r_{21}R_{41}\bar{z}_{31}+r_{31}R_{41}\bar{z}_{21}\\
r_{21}\bar{z}_{31}\bar{z}_{41}+r_{31}\bar{z}_{21}\bar{z}_{41}
+R_{41}\bar{z}_{21}\bar{z}_{31}\\
\bar{z}_{21}\bar{z}_{31}\bar{z}_{41}\end{array}\!\right\rgroup\]
\[v_2=\frac{1}{\sqrt{r_{32}R_{42}r_{21}}}\left\lgroup\!\begin{array}{c}
-r_{32}R_{42}z_{21}\\-r_{32}z_{21}\bar{z}_{42}+r_{32}R_{42}r_{21}
-z_{21}R_{42}\bar{z}_{32}\\
r_{32}r_{21}\bar{z}_{42}-z_{21}\bar{z}_{32}\bar{z}_{42}
+R_{42}\bar{z}_{32}r_{21}\\
\bar{z}_{32}r_{21}\bar{z}_{42}\end{array}\!\right\rgroup\]
\[v_3=\frac{1}{\sqrt{R_{43}r_{31}r_{32}}}\left\lgroup\!\begin{array}{c}
z_{31}z_{32}R_{43}\\
z_{31}z_{32}\bar{z}_{43}-z_{31}R_{43}r_{32}-z_{32}R_{43}r_{31}\\
-z_{31}r_{32}\bar{z}_{43}-z_{32}r_{31}\bar{z}_{43}+R_{43}r_{31}r_{32}\\
r_{31}r_{32}\bar{z}_{43}\end{array}\!\right\rgroup\]
\[v_4=\frac{1}{\sqrt{R_{41}R_{42}R_{43}}}\left\lgroup\!
\begin{array}{c}
-z_{41}z_{42}z_{43}\\
z_{41}z_{42}R_{43}+z_{41}z_{43}R_{42}+z_{42}z_{43}R_{41}\\
-z_{41}R_{42}R_{43}-z_{42}R_{41}R_{43}-z_{43}R_{41}R_{42}\\
R_{41}R_{42}R_{43}\end{array}\!\right\rgroup\]
and we may take $M$ to be the matrix with column vectors~$v_i$.  Hence,
\[\det M=P/(r_{21}r_{31}r_{32}R_{41}R_{42}R_{43})\]
where $P$ is a polynomial consisting of monomials each of which
contains one of $r_{ij}^2$, $r_{ij}z_{ij}$, $r_{ij}\bar{z}_{ij}$, or
$z_{ij}\bar{z}_{ij}$ for each $j<i<4$, and one of $R_{4j}^2$,
$R_{4j}z_{4j}$, $R_{4j}\bar{z}_{4j}$, or $z_{4j}\bar{z}_{4j}$.  Since
$z_{ij}\bar{z}_{ij}=r_{ij}^2$, each monomial is divisible by $r_{ij}$
and, since $|z_{4j}|^2=R_{4j}(r_{4j}-r)$, each monomial is divisible
by $R_{4j}$. Therefore, we can divide by the factor of
$r_{21}r_{31}r_{32}R_{41}R_{42}R_{43}$ leaving monomials with one of
$r_{ij}$, $z_{ij}$ or $\bar{z}_{ij}$ for each $j<i<4$, and one of
$(r_{4j}+r)$, $z_{4j}$, $\bar{z}_{4j}$ or $(r_{4j}-r)$.  It follows that
$\det M$ is now expressed as a homogeneous degree 6 polynomial in $r$,
$r_{ij}$, $z_{ij}$ and $\bar{z}_{ij}$ for $j<i\leq 4$.

Recall that $\det M$ is invariant under Euclidean motions.  Moreover, the six
distances $r_{ij}$ determine our configuration of four points.  Also, notice
that these distances are constrained only by triangle inequalities.  Hence, the
Atiyah determinant $\det M$ may be regarded as a function of the independent
variables~$r_{ij}$.  We claim that
\begin{itemize}
\item $\Re(\det M)$ is a polynomial in~$r_{ij}$ (homogeneous of degree 6).
\item $|\det M|^2$ is a polynomial in~$r_{ij}$ (homogeneous of degree 12).
\end{itemize}
Having done this we shall use the triangle inequalities on the four
faces of our configuration to show that, in fact, $\Re(\det M)>0$.
This is more than sufficient to finish the proof.

It is convenient to set $z_{ij}=-z_{ji}$ and $r_{ij}=r_{ji}$ when
$j>i$.  The monomials in our expression for $\det M$ contain an equal
number of $z_{ij}$ and $\bar{z}_{kl}$.  Consider the product
$z_{ij}\bar{z}_{kl}$.  There are two cases: \newline
\indent
(i) $\{i,j\}$ and $\{k,l\}$ have an element in common (suppose $l=j$):
\begin{equation}  \label{eq:adef}
    z_{ij}\bar{z}_{kj}=(1/2)\left(r_{ij}^2+r_{kj}^2-r_{ki}^2\right)+
    2A_{ijk}\sqrt{-1}
\end{equation}
where $A_{ijk}$, defined by~(\ref{eq:adef}), equals plus or minus
the area of the $ijk^{\mathrm{th}}$ triangle under the projection
${\mathbb R}^3={\mathbb R}\times{\mathbb C}\to{\mathbb C}$ onto the
complex plane; \newline
\indent(ii) $\{ i,j,k,l\}=\{1,2,3,4\}$:
\[z_{ij}\bar{z}_{kl}=z_{ij}\bar{z}_{kj}z_{kj}\bar{z}_{kl}/r_{kj}^2\]
and we can rewrite the numerator as in (i).\newline
We claim that all quadratic expressions in the $A_{ijk}$ may be
written as polynomials in the $r_{ij}$ and $r^2$.  Specifically, when
all four points lie in the complex plane, then one may verify that
\begin{equation}\label{quad}
16A_{ijk}A_{ijl}=2r_{ij}^2(r_{ik}^2+r_{il}^2-r_{kl}^2)
-(r_{ij}^2+r_{ik}^2-r_{jk}^2)(r_{ij}^2+r_{il}^2-r_{jl}^2)
\end{equation}
and when the fourth point lies off the plane ($r>0$), we replace
$r_{4j}^2$ by $r_{4j}^2-r^2$.

Now, observe that our formulae so far for $\Re(\det M)$ and $(\Im(\det M))^2$
involve only quadratics expressions in $A_{ijk}$.  If we substitute according
to (\ref{quad}) and its non-planar version, we obtain rational expressions for
$\Re(\det M)$ and $(\Im(\det M))^2$ in the seven quantities in $r_{ij}$
and $r$, the denominator being a polynomial in the six
variables~$r_{ij}$.  Recall that a reflection such as $r\mapsto -r$ conjugates
$\det M$.  Hence, we may drop all odd powers of $r$ in the numerators, to
obtain polynomials in $r_{ij}$ and~$r^2$.

Finally, we eliminate~$r^2$ from these expressions.  This is possible
by writing the volume $V$ of the tetrahedron with vertices
our four points in two different ways.  On the one hand
\begin{eqnarray*}
    144V^2&=&-r_{21}^4r_{43}^2-r_{21}^2r_{43}^4-r_{32}^2r_{41}^4
    -r_{32}^4r_{41}^2-r_{31}^4r_{42}^2-r_{31}^2r_{42}^4\\
    &&+r_{21}^2r_{43}^2r_{31}^2+r_{21}^2r_{43}^2r_{41}^2
   -r_{21}^2r_{42}^2r_{41}^2+r_{21}^2r_{42}^2r_{43}^2\\
   &&+r_{21}^2r_{42}^2r_{31}^2+r_{21}^2r_{32}^2r_{43}^2
   -r_{21}^2r_{32}^2r_{31}^2+r_{32}^2r_{42}^2r_{41}^2\\
   &&+r_{31}^2r_{42}^2r_{41}^2+r_{32}^2r_{43}^2r_{41}^2
   -r_{32}^2r_{42}^2r_{43}^2+r_{32}^2r_{42}^2r_{31}^2\\
   &&+r_{31}^2r_{42}^2r_{43}^2+r_{32}^2r_{31}^2r_{41}^2
   -r_{31}^2r_{43}^2r_{41}^2+r_{21}^2r_{32}^2r_{41}^2.
\end{eqnarray*}
On the other hand, let $A$ denote the area of the
triangular base in~${\mathbb C}$. Then
\[ 16A^2=2r_{21}^2r_{31}^2+2r_{21}^2r_{32}^2+2r_{31}^2r_{32}^2-r_{21}^4
-r_{31}^4-r_{32}^4\]
and $V=rA/3$.  We may therefore replace $r^2$ by
$9V^2/A^2$, as required.

Thus, we may conclude that $\Re(\det M)$ and $(\Im(\det M))^2$ are rational
functions of the variables $r_{ij}$ but, if we now clear any common factors, we
claim they are, in fact, polynomials.  To see this, notice that in (ii) we had
a choice when we introduced $\bar{z}_{kj}z_{kj}/r_{kj}^2$.  We could have
insisted that $\{k,j\}\subset\{1,2,3\}$.  Then the denominators of
$\Re(\det M)$ and $(\Im(\det M))^2$ would not involve~$r_{4j}$.  However,
$\det M$ does not see the ordering of our four points.  Hence, if $r_{4j}$ are
omitted from the denominators, then so are all variables~$r_{ij}$, as required.

It remains to calculate these polynomials.  We did this using
Maple\footnote{The program is at:
{\footnotesize\tt ftp://ftp.maths.adelaide.edu.au/meastwood/maple/points}} and found
that $\Re(\det M)$ is a homogeneous polynomial of degree 6 with 226
terms and $|\det M|^2$ is a homogeneous degree 12 polynomial with 4500
terms.  We claim that $\Re(\det M)>0$.  To see this, we can rewrite
the output of the Maple calculation as follows:
\begin{eqnarray*}
\Re(\det M)&\!=\!&64r_{21}r_{31}r_{32}r_{41}r_{42}r_{43}
    -4d_3(r_{21}r_{43},r_{31}r_{42},r_{32}r_{41})\\
    &&\;+12\av\left(r_{41}((r_{42}+r_{43})^2-r_{32}^2)
    d_3(r_{21},r_{31},r_{32})\right)+288V^2.
\end{eqnarray*}
Here, $d_3(a,b,c)$ is the polynomial~(\ref{deethree}) and $\av$
denotes the operation of averaging a polynomial in $r_{ij}$ under the action
of ${\mathcal S}_4$ on the vertices of our tetrahedron: for example,
$$\begin{array}{rcl}
\av(r_{21})&=&(r_{21}+r_{31}+r_{32}+r_{41}+r_{42}+r_{43})/6\\[4pt]
\av(r_{21}r_{43})&=&(r_{21}r_{43}+r_{31}r_{42}+r_{41}r_{32})/3.\end{array}$$
The final two terms are non-negative since the triangle inequality gives
$$(r_{42}+r_{43})^2\geq r_{32}^2\quad\mbox{and}\quad
d_3(r_{21},r_{31},r_{32})\geq 0,$$
and the square of the volume is non-negative.  To estimate the other terms we
may use the easily verified inequality
\[ abc\geq d_3(a,b,c),\quad \forall\ a,b,c\geq 0.\]
In conclusion,
\[\Re(\det M)\geq 60r_{21}r_{31}r_{32}r_{41}r_{42}r_{43}>0,\]
as required.  This is nearly enough for a stronger conjecture of
Atiyah and Sutcliffe \cite[Conjecture 2]{as} that
$|\det M|\geq 64r_{21}r_{31}r_{32}r_{41}r_{42}r_{43}$.
\end{proof}
A third conjecture of Atiyah and Sutcliffe \cite[Conjecture 3]{as} can
be expressed in the four point case in terms of polynomials in the
edge lengths as:
\[|\det M|^2\geq \prod_1^4 (d_3(r_{ij},r_{ik},r_{jk})+8r_{ij}r_{ik}r_{jk})\]
where the product runs over the four faces of the tetrahedron and the
left hand side is known explicitly.  We have been unable to prove this
conjecture even in the case that the four points lie in a plane (in which
case $|\det M|$ can be replaced by the simpler expression $\Re(\det M)$
given above).

\section{The planar case}
Atiyah's basic conjecture is unresolved in general,
even when the $n$ points lie in a
plane (in which case recall that $\det M$ is real).  Reasoning analogous to the
case of four points gives the following.
\begin{thm}
    The Atiyah determinant of $n$ points in a plane can be expressed
    as a rational function in the distances between the points.
\end{thm}
\proof
    Again, we can express $z_{ij}\bar{z}_{kl}$ as a rational function in
    the $r_{ij}$ and $A_{ijk}$.  It is no longer true that quadratic
    expressions in the $A_{ijk}$ are polynomials in the $r_{ij}$.
    Instead they are rational functions in the $r_{ij}$.  This uses the
    same trick of introducing new points in common between two
    triangles in order to apply~(\ref{quad}):
    $$ A_{ijk}A_{lmn}=\frac{(A_{ijk}A_{ijn})(A_{imn}A_{lmn})}
    {A_{ijn}A_{imn}}.\eqno{\qed}$$

In the general four point case, the distances $r_{ij}$ acted as variables.  The
denominator was too small to be appropriately symmetrical and therefore had to
divide the numerator, leaving a polynomial rather than a rational function.  In
the planar case (and also in the general case with more than four points), the
distances between points satisfy a set of polynomial constraints.
Symmetry arguments are no longer valid and expressions for the determinant are
no longer unique.  We suspect, however, that there is a polynomial expression.

\rk{Acknowledgements}Support from the Australian Research Council and
the Mathematical Sciences Research Institute is gratefully
acknowledged.  Research at MSRI is supported in part by NSF grant
DMS-9701755.

\end{document}